\documentclass{amsart}

\RequirePackage{amsmath,amssymb}

\RequirePackage{xypic}
\xyoption{all}
\CompileMatrices

\theoremstyle{plain}

\newtheorem{theorem}{Theorem}[section]
\newtheorem{proposition}[theorem]{Proposition}
\newtheorem{lemma}[theorem]{Lemma}
\newtheorem{corollary}[theorem]{Corollary}
\newtheorem{conjecture}[theorem]{Conjecture}
\newtheorem{mthm}{Theorem}

\theoremstyle{definition}

\newtheorem{defn}[theorem]{Definition}
\newtheorem{example}[theorem]{Example}

\renewcommand{\H}{\mathcal{H}}
\renewcommand{\O}{\mathcal{O}}
\newcommand{\Q}{\mathbf{Q}}
\renewcommand{\S}{\mathcal{S}}
\newcommand{\Z}{\mathbf{Z}}

\newcommand{\Fi}{F_{\infty}}
\newcommand{\Fiw}{F_{\infty,w}}
\newcommand{\Fn}{F_{n}}
\newcommand{\Fnw}{F_{n,w}}
\newcommand{\FS}{F_{\Sigma}}
\newcommand{\Fv}{F_{v}}
\newcommand{\Fw}{F_{w}}
\newcommand{\Qi}{\Q_{\infty}}
\newcommand{\Qiw}{\Q_{\infty,w}}
\newcommand{\Qp}{\Q_{p}}
\newcommand{\Zp}{\Z_{p}}

\newcommand{\GF}{G_{F}}
\newcommand{\GFi}{G_{\Fi}}
\newcommand{\Gl}{G_{\ell}}

\newcommand{\GQ}{G_{\Q}}
\newcommand{\Gv}{G_{v}}
\newcommand{\Gw}{G_{w}}
\newcommand{\Il}{I_{\ell}}
\newcommand{\Iv}{I_{v}}
\newcommand{\Iw}{I_{w}}

\DeclareMathOperator{\crit}{cr}
\newcommand{\cyc}{\varepsilon}

\newcommand{\dvm}{d_{v}^{-}}

\newcommand{\HfS}{H_{f,\S}^{1}}
\newcommand{\HfSA}{H_{f,\S(A)}^{1}}
\newcommand{\HfSm}{H_{f,\Sm}^{1}}
\newcommand{\Hs}{H_{s}^{1}}
\newcommand{\Hsc}{H_{s,\crit}^{1}}
\newcommand{\HsS}{H_{s,\S}^{1}}

\DeclareMathOperator{\Sel}{Sel}
\newcommand{\Sm}{\S_{\min}}
\newcommand{\Selc}{\Sel_{\crit}}
\newcommand{\SelS}{\Sel_{\S}}
\newcommand{\SelSm}{\Sel_{\Sm}}
\newcommand{\SelSA}{\Sel_{\S(A)}}

\newcommand{\Ab}{\Delta}

\newcommand{\rhob}{\bar{\rho}}

\newcommand{\inj}{\hookrightarrow}

\newcommand{\too}{\longrightarrow}

\newcommand{\Ar}{A_{\rho}}
\newcommand{\At}{A_{t}}
\newcommand{\Atw}{A_{t,w}}

\newcommand{\Tr}{T_{\rho}}
\newcommand{\Vr}{V_{\rho}}
\newcommand{\Vrv}{V_{\rho,v}}

\DeclareMathOperator{\ad}{ad}
\DeclareMathOperator{\B}{B}
\DeclareMathOperator{\G}{G}
\DeclareMathOperator{\GL}{GL}
\DeclareMathOperator{\GSp}{GSp}
\DeclareMathOperator{\mult}{mult}

\DeclareMathOperator{\ZG}{Z}

\newcommand{\Blie}{\mathcal B}
\newcommand{\Ghat}{\widehat{\G}}
\newcommand{\Glie}{\mathcal G}

\DeclareMathOperator{\coker}{coker}
\DeclareMathOperator{\corank}{corank}
\DeclareMathOperator{\Gal}{Gal}
\DeclareMathOperator{\Hom}{Hom}
\DeclareMathOperator{\im}{im}
\DeclareMathOperator{\no}{no}
\DeclareMathOperator{\tors}{tors}
\DeclareMathOperator{\tr}{tr}

\begin{document}

\title{Iwasawa invariants of Galois deformations}
\author{Tom Weston}
\email{taweston@amherst.edu}
\address{Department of Mathematics and Statistics, University of Massachusetts Amherst}

\maketitle

Let $p$ be an odd prime and let $K$ be a finite extension of $\Qp$
with residue field $k$.
Let $\G$ be a reductive group over the ring of integers of $K$ and fix
a continuous representation
$$\rhob : \GF \to \G(k)$$
of the absolute Galois group of a number field $F$.
Assume that $\rhob$ is ordinary in the sense that the image of
any decomposition group at a place $v$ dividing $p$ lies in
some Borel subgroup $\B_{v}$ of $\G$.
Assume also that $\rhob$ satisfies the conditions of
\cite[Section 7]{Tilouine} which guarantee that it has a reasonable deformation
theory; see Section~\ref{s31} for details.
In this paper we show that the Iwasawa invariants of the Selmer
group of a nearly ordinary
deformation of $\rhob$
depends only on $\rhob$ and the tame ramification of the deformation.

For a precise statement,
let $\H$ denote the set of nearly ordinary (with respect to the $\B_{v}$)
deformations of $\rhob$ to continuous finitely ramified representations
$$\rho : \GF \to \G(\O^{\tr})$$
over the ring of integers $\O^{\tr}$ of the maximal totally ramified
extension of $K$.  Fix an algebraic representation $r : \G \to \GL_{n}$
such that $[F:\Q]$ divides
$$\sum_{v\text{~real}} \dim (r \circ \rhob)^{c_{v}=-1}
+ \sum_{v\text{~complex}} n$$
(where $c_{v}$ is complex conjugation at $v$)
and such that $r \circ \rhob$ and its Cartier dual have trivial
$\GF$-invariants.
Subject to the choice of Borel subgroups $\tilde{\B}_{v}$ of $\GL_{n}$
containing $r(\B_{v})$,
for any $\rho \in \H$ we define a Selmer group
$\Sel(\Fi,\rho,r)$ over the cyclotomic $\Zp$-extension $\Fi$ of $F$.  
This Selmer group is closely related to the usual Selmer group of
Greenberg and thus is conjecturally related to the $p$-adic
$L$-function of $\rho$ with respect to $r$.

We write $\H(r)$ for the set of $\rho \in \H$ for which
$\Sel(\Fi,\rho,r)$ is cotorsion over
the Iwasawa algebra $\O^{\tr}[[\Gal(\Fi/F)]]$ and which satisfy a certain
ramification condition; see Section~\ref{s32} for
a precise definition.  For $\rho \in \H(r)$ we say that
$\mu(\rho,r)$ vanishes if $\Sel(\Fi,\rho,r)$ is cofinitely generated over
$\O^{\tr}$ and we let $\lambda(\rho,r)$ denote the $\O^{\tr}$-corank
of $\Sel(\Fi,\rho,r)$.
Our main result is the following theorem on the variation of
these Iwasawa invariants over $\H(r)$.  For a place $v$ of $F$ set
$$\delta_{v}(\rho,r) = \sum_{w \mid v}
\mult_{\omega}(r \circ \rhob|_{\Iw}) -
\mult_{\omega}\bigl((r \circ \rho) \otimes K|_{\Iw}\bigr);$$
here the sum runs over the places of $\Fi$ dividing $v$ and
$\mult_{\omega}(\cdot)$ is the multiplicity of the
Teichm\"uller character in the given representation of $\Gw/\Iw$.
Note that $\delta_{v}(\rho,r)$ depends only on $r \circ \rhob$ and
the restriction of $r \circ \rho$ to $\Iv$.

\begin{mthm} \label{mt1}
If $\mu(\rho_{0},r) = 0$ for some $\rho_{0} \in \H(r)$, then
$\mu(\rho,r) = 0$ for all $\rho \in \H(r)$.  If this is the case,
then the difference
\begin{equation} \label{eq:diff}
\lambda(\rho,r) - \sum_{v \nmid p} \delta_{v}(\rho,r)
\end{equation}
is independent of $\rho \in \H(r)$.
\end{mthm}

In particular, Theorem~\ref{mt1} implies that the values 
$\lambda(\rho,r)$ for all $\rho \in \H(r)$
can be easily determined from the knowledge of 
$\lambda(\rho_{0},r)$ for a single $\rho_{0}$ with $\mu(\rho_{0},r)=0$.
The question of which $\lambda$-invariants can
occur in $\H(r)$
is intimately related to the questions of level raising and lowering for
$\rhob$.  
For example, if level lowering holds for $\rhob$, 
then (\ref{eq:diff}) equals the minimal
$\lambda$-invariant in the family $\H(r)$.

Of course, although conjecturally $\H(r)$ is quite large, this is known
in very few cases.  The case of $F = \Q$, $\G = \GL_{2}$ and $r$
the identity is studied via Hida theory in \cite{EPW}.
The work of \cite{fujiwara} and
\cite{hida3} provide additional cases where we can
apply our results.  Specifically, let $F$ be totally real and let
$\G = \GL_{2}$.  
Let $\rhob : \GF \to \GL_{2}(k)$ be an absolutely irreducible
residual representation attached to
some $p$-ordinary Hilbert modular form $f_{0}$ over $F$ which is not
of CM-type.  Assume that:
\begin{itemize}
\item $\rhob$ is absolutely irreducible when restricted to
$F(\sqrt{(-1)^{(p-1)/2}p})$;
\item $F \cap \Q(\mu_{p}) = \Q$;
\item $(\O_{F} \otimes_{\Z} \Zp)^{\times}$ has no $p$-torsion;
\item $F/\Q$ is unramified at $p$.
\end{itemize}
If $\ad^{0} : \GL_{2} \to \GL_{3}$ is the trace zero adjoint representation,
then $\H(\ad^{0})$ contains the Galois representations associated to the set
$\H^{0}$ of $p$-ordinary
Hilbert modular forms $f$ which are congruent to $f_{0}$ in $k$.
In this case our critical Selmer group coincides with Greenberg's ordinary
Selmer group, so that we obtain the following result on the classical
Iwasawa invariants of forms in $\H^{0}$.

\begin{mthm} \label{mt2}
If $\mu(\ad^{0}f_{0})=0$,
then $\mu(\ad^{0}f)=0$ for all $f \in \H^{0}$.
If this is the case, then the difference
$$\lambda(\ad^{0}f) - \sum_{v \nmid p} 
\delta_{v}(\ad^{0}f)$$
is independent of $f \in \H^{0}$.
\end{mthm}

The origin of this work is the paper \cite{GreenbergVatsal}
of Greenberg and Vatsal, in which analogous results were obtained for
Tate modules of elliptic curves.  These results were extended to
arbitrary modular forms over $\Q$ in \cite{EPW}.
The methods given in this paper are a generalization of those
of \cite{EPW} to arbitrary algebraic groups and number fields.

The first issue we must confront for this generalization is that
variation of Hodge--Tate weights causes
Greenberg's Selmer groups of ordinary representations
to behave poorly in families.  In Section~\ref{s1}
we introduce and study
our critical Selmer groups, which are essentially Greenberg's
Selmer group with extra twists to facilitate interpolation over families.

Distinct 
deformations of a fixed residual representation may have very different
Selmer groups.  However, in favorable situations it is possible to
visualize these differences in the cohomology of the residual
representation.  This is done via the theory of residual Selmer groups
developed in Section~\ref{s2}.  We also discuss the connections with
level lowering and level raising, focusing on the less conjectural
case of $\GL_{2}$.
We obtain Theorem~\ref{mt1}
in Section~\ref{s3} by applying the preceding results to the study
of Iwasawa invariants of nearly ordinary
deformations of fixed residual representations as in \cite{Tilouine}.

This work grew out of the paper \cite{EPW};
it is a pleasure to thank Matthew Emerton and Robert Pollack
for many helpful discussions on this material.  The author would also
like to thank Paul Gunnells, Chris Skinner and Eric Sommers 
for their help with this project.

\section{Selmer groups} \label{s1}

Throughout this paper $p$ denotes an odd prime.

\subsection{Local Galois representations} \label{s11}

Let $K$ and $L$ denote finite extensions of $\Qp$ and
let $G_{L}$ denote the absolute Galois group of $L$.
A {\it nearly ordinary $G_{L}$-representation} over $K$ is a
finite-dimensional $K$-vector space $V$ endowed with
a continuous $K$-linear action of $G_{L}$ and a choice of a
$G_{L}$-stable complete flag
$$0 = V^{0} \subsetneq V^{1} \subsetneq \cdots \subsetneq
V^{n} = V.$$
Let $\chi_{i} : G_{L} \to K^{\times}$ denote the character by 
which $G_{L}$ acts on $V^{i}/V^{i-1}$.
If $V$ is a Hodge--Tate representation of $G_{L}$, then $\chi_{i}$ must be the
product of a character of finite order and 
some integer power $\cyc^{m_{i}}$ of the cyclotomic character.
In this case we 
call $m_{1},\ldots,m_{n}$ the {\it Hodge--Tate weights} of $V$.

We say that a nearly ordinary $G_{L}$-representation $V$ is {\it ordinary} if
it is Hodge--Tate and if
$$m_{1} \geq m_{2} \geq \cdots \geq m_{n}.$$
(Note that as we have defined it the property of being ordinary depends
on the choice of complete flag and not merely on $V$.)
It is well known (see \cite{PerrinRiou}) that ordinary representations
are always potentially semistable.
In fact, the converse is essentially true as well.

\begin{lemma}
If $V$ is nearly ordinary and potentially semistable, then there
exists a complete $G_{L}$-stable flag
$$0 = \tilde{V}^{0} \subsetneq \tilde{V}^{1} \subsetneq \cdots \subsetneq
\tilde{V}^{n} = V$$
with respect to which $V$ is ordinary.
\end{lemma}
\begin{proof}
Let
$$0 = V^{0} \subsetneq V^{1} \subsetneq \cdots \subsetneq
V^{n} = V$$
be the given $G_{L}$-stable complete flag.  $V$ is Hodge--Tate since it
is potentially semistable;
let $m_{1},\ldots,m_{n} \in \Z$ denote the Hodge--Tate
weights as before.  Let $i$ be the least index
such that $m_{i} < m_{i+1}$ and let $W$ denote the two-dimensional
$G_{L}$-representation $V^{i+1}/V^{i-1}$.  Choosing a basis $x,y$ for
$W$ with $x \in V^{i}/V^{i-1}$, the representation of $G_{L}$ on $W$ has
the form
$$\left( \begin{array}{cc} \cyc^{m_{i}} & * \\ 0 & \cyc^{m_{i+1}}
\end{array} \right)$$
and is potentially semistable (as this property is preserved under
passage to subquotients).  Thus the $*$ above may be regarded
as an element of $H^{1}_{g}(\Qp,K(m_{i}-m_{i+1}))$, in the notation
of \cite{BlochKato}.  This group is trivial since $m_{i} - m_{i+1} < 0$,
so that in fact the representation of $G_{L}$ on
$V^{i+1}/V^{i-1}$ is diagonal.  It follows that the complete flag
$$\tilde{V}^{j} = \begin{cases} V^{j} & j < i \\
V^{i-1} + K \cdot y & j = i \\
V^{i+1} & j= i+1
\end{cases}$$
of $V^{i+1}$ gives $V^{i+1}$ the structure of ordinary
$G_{L}$-representation.  Continuing in this way yields the desired
ordinary flag for $V$.
\end{proof}

\subsection{Global Galois representations} \label{s12}

Let $F$ be a finite extension of $\Q$ and let $\Fi$ denote the
cyclotomic $\Zp$-extension of $F$.
A {\it nearly ordinary Galois representation} over a finite extension
$K$ of $\Qp$ is a finite-dimensional
$K$-vector space equipped with a $K$-linear action of the
absolute Galois group $\GF$ such that
$V$ is equipped with the structure
$$0 = V^{0}_{v} \subsetneq V^{1}_{v} \subsetneq \cdots \subsetneq
V^{n}_{v} = V$$
of nearly ordinary $\Gv$-representation
for each place $v$ of $F$ dividing $p$.  A nearly ordinary Galois
representation is said to be {\it ordinary} if it is ordinary at
each place $v$ dividing $p$ and if the corresponding
Hodge--Tate weights for each $v$ coincide.

Let $V$ be a nearly ordinary Galois representation of dimension $n$.
For a real place $v$ of $F$, let $\dvm(V)$
denote the $K$-dimension of the subspace of $V$ on which complex
conjugation at $v$ acts by
$-1$.  We say that $V$ is {\it critical} if
$$\sum_{v \text{~real}} \dvm(V) + \sum_{v \text{~complex}} n$$
is divisible by $[F:\Q]$.  (Note that this is automatic if $F=\Q$ or
if $F$ is totally complex and $n$ is even.)
If $c(V)$ denotes this quotient, we then set
$$V_{v}^{\crit} = V_{v}^{n-c(V)}$$
for each place $v$ dividing $p$.

\begin{example}
Let $F$ be totally real and let $V$ be the two-dimensional representation
associated to a $p$-ordinary Hilbert modular form $f$ of parallel weight
$(k,k,\ldots,k)$ with $k \geq 2$.  Then $V$ has a natural structure of
ordinary Galois representation with Hodge--Tate weights $0$ and $k-1$.
We have $\dvm(V) = 1$ for every archimedean $v$, so that $V$ is
critical with $c(V)=1$.  If $\ad^{0}V$ denotes the trace zero adjoint of $V$,
then $\ad^{0}V$ is ordinary with Hodge--Tate weights $1-k,0,k-1$.
We have $\dvm(\ad^{0}V)=2$ for each $v$, so that $V$ is critical with
$c(V)=2$.
\end{example}

\begin{example}
Let $\rho : \GQ \to \GL_{n}(K)$ be associated to a cohomological
cuspidal representation of $\GL_{n}$ of highest weight
$(a_{1},a_{2},\ldots,a_{n})$ with $a_{1} \geq a_{2} \geq \cdots \geq
a_{n}$ and let $V$ denote
the associated $n$-dimensional Galois representation.
If $\rho$ is ordinary, then $V$ has a natural structure of ordinary
Galois representation; the Hodge--Tate weights are conjectured to be
$$a_{1}+n-1,a_{2}+n-2,a_{3}+n-3,\ldots,a_{n}.$$
\end{example}

\begin{example}
Let $\rho : \GQ \to \GSp_{4}(K)$ be associated to a $p$-ordinary cohomological
cuspidal representation of $\GSp_{4}$ of highest weight
$(a,b;a+b)$ with $a \geq b \geq 0$ 
as in \cite{TilouineUrban} and let $V$ denote
the associated four-dimensional Galois representation; it 
has a natural structure of ordinary Galois representation
with Hodge--Tate weights
$$a+b+3,a+2,b+1,0.$$
We have
$\dvm(V)=2$ for the unique archimedean place of $\Q$, so that
$c(V)=2$.
\end{example}

\subsection{Critical Selmer groups} \label{s13}

Let $V$ be a critical nearly 
ordinary Galois representation as above.
Let $\O$ denote the ring of integers of $K$.
Fix a $\GF$-stable $\O$-lattice $T$ in $V$ and set $A = V/T$; we call
$A$ a {\it torsion quotient} of $V$, although in general $A$
is not uniquely determined by $V$.
We say that a finite set of places $\Sigma$ of $F$
is {\it sufficiently large} for $A$ if it contains all archimedean places,
all places dividing $p$, and all places at which $A$ is ramified.

For any place $w$ of $\Fi$ set
$$\Hsc(\Fiw,A) = \begin{cases} H^{1}(\Gw,A) & w \nmid p; \\
\im H^{1}(\Gw,A) \to H^{1}(\Iw,A/A_{w}^{\crit}) & w \mid p; \end{cases}$$
where $A_{w}^{\crit}$ is the image of $V_{v}^{\crit}$ in $A$ with $v$
the restriction of $w$ to $F$.
We define the {\it critical Selmer group} of $A$ by
\begin{align}
\Selc(\Fi,A) &= \ker H^{1}(\Fi,A) \to \prod_{w} \Hsc(\Fiw,A) \notag \\
&= \ker H^{1}(\FS/\Fi,A) \to \prod_{w \mid v 
\in \Sigma} \Hsc(\Fiw,A) \label{eq:def}
\end{align}
for any finite set of places $\Sigma$ of $F$
which is sufficiently large for $A$.

By \cite[Propositions 1 and 2]{Greenberg1} we have
$$\corank_{\Lambda}\Hsc(\Fiw,A) =
\begin{cases} 0 & w \nmid p; \\
[F_{v}:\Qp]c(V) & w \mid v \mid p.\end{cases}$$
By \cite[Proposition 3]{Greenberg1} we also have
\begin{multline*}
\corank_{\Lambda} H^{1}(\FS/\Fi,A) =
\sum_{v \text{~real}} \dvm(V) + \sum_{v \text{~complex}} n + \corank_{\Lambda}
 H^{2}(\FS/\Fi,A) \\
=[F:\Q] \cdot c(V) + \corank_{\Lambda} H^{2}(\FS/\Fi,A)
\end{multline*}
for any finite set of places $\Sigma$ sufficiently large for $A$;
here $\Lambda = \O[[\Gal(\Fi/F)]]$ is the Iwasawa algebra.
Since $\sum_{v\mid p} [F_{v}:\Qp] = [F:\Q]$,
by (\ref{eq:def}) 
(together with \cite[Proposition 6]{Greenberg1})
we thus have the following result.

\begin{proposition} \label{prop:small}
$\Selc(\Fi,A)$ is a cofinitely generated $\Lambda$-module of
$\Lambda$-corank at least $\corank_{\Lambda} H^{2}(\FS/\Fi,A)$.
\end{proposition}

In fact, it is conjectured \cite{Greenberg1} 
that the error term above vanishes.

\begin{conjecture}[Greenberg]
For $A$ as above $H^{2}(\FS/\Fi,A)=0$.
\end{conjecture}

It thus seems reasonable to adapt
standard conjectures on Selmer groups as follows.

\begin{conjecture}
The critical Selmer group
$\Selc(\Fi,A)$ is $\Lambda$-cotorsion.
\end{conjecture}

In fact, we will see below that the latter conjecture implies the
former.

If $V$ is ordinary with 
distinct Hodge--Tate weights $m_{n-c(V)}$ and $m_{n-c(V)+1}$,
then the critical Selmer group is closely related to
Greenberg's Selmer group.  More precisely, for any $m$ such that
$$m_{n-c(V)+1} \leq m < m_{n-c(V)},$$
we have
$$\Selc(\Fi,A) = \Sel(\Fi,A(m) \otimes \omega^{-m})$$
where the latter Selmer group is as in
\cite[Section 7]{Greenberg1},
$A(m) = A \otimes \cyc^{m}$ is the twist of $A$ by $m$ powers
of the cyclotomic character and $\omega$ is the Teichm\"uller character.
Note that $A(m) \otimes \omega^{-m} \cong A$
as $\GFi$-modules; the twist above is simply shifting the Hodge
filtration.

\subsection{Structure of critical Selmer groups} \label{s14}

Let $A$ be a torsion quotient of a critical nearly ordinary Galois
representation $V$ as before.
Greenberg \cite{Greenberg2, GreenbergVatsal} has obtained
powerful results on the structure of
Selmer groups under the assumption that they are $\Lambda$-cotorsion.
In this section we adapt his techniques to critical Selmer groups.
The essential idea is to replace $A$ by a twist which is well-behaved.
Let 
$$\kappa = \omega^{-1}\cyc : \GF \to 1 + p\Zp$$ 
denote the character
giving the isomorphism $\Gal(\Fi/F) \cong 1+p\Zp$.  For any $t \in \Zp$
we set $\At = A \otimes \kappa^{t}$ and
$A_{t,w}^{\crit} = 
A^{\crit}_{w} \otimes
\kappa^{t}$ for any place $w$ of $\Fi$ dividing $p$.
Note that $\At$ is
isomorphic to $A$ as a $\GFi$-module.
Let $\At^{*}$ denote the Cartier dual
$\Hom_{\O}(\At,K/\O(1))$; it is
a free $\O$-module of finite rank.

\begin{proposition} \label{prop:loc}
Let $A$ be a torsion quotient of a
critical nearly ordinary Galois representation $V$.  Assume that
$\Selc(\Fi,A)$ is $\Lambda$-cotorsion and that $H^{0}(\Fi,A^{*} \otimes K/\O)$
is finite.
\begin{enumerate}
\item The sequence
$$0 \to \Selc(\Fi,A) \to H^{1}(\FS/\Fi,A) \to \prod_{w \mid 
v \in \Sigma} \Hsc(\Fiw,A)
\to 0$$
is exact for any finite set of places $\Sigma$ of $F$
sufficiently large for $A$.
\item Assume that $H^{0}(F,A^{*} \otimes K/\O)=0$.  Then $\Selc(\Fi,A)$ has no
proper $\Lambda$-submodules of finite index.
\end{enumerate}
\end{proposition}
\begin{proof}
The sequence in (1) is exact by definition except for the surjectivity on
the right.  In fact $\Hsc(\Fiw,A)$ is divisible for each $w$ since
$\Gw$ has $p$-cohomological dimension one, so
to prove (1) it suffices to show that the the cokernel of the map
$$\gamma : H^{1}(\FS/\Fi,A) \to \prod_{w \mid 
v \in \Sigma} \Hsc(\Fiw,A)$$
is finite.

Let $\Fn$ denote the unique subfield of $\Fi$ of degree $p^{n}$ over
$F$.  For a place $w$ of $\Fn$ and $t \in \Zp$ set
$$\Hs(\Fnw,\At) = \begin{cases} H^{1}(\Iw,\At)^{\Gw/\Iw} & w \nmid p \\
\im H^{1}(\Gw,\At) \to H^{1}(\Gw,\At/A_{t,w}^{\crit}) & w \mid p. \end{cases}$$
(The condition at places dividing $p$ is slightly stronger
than we usually impose so that we may apply
the appropriate duality results below.)
We claim that 
to prove the surjectivity of $\gamma$
it suffices to show that for some $t$ the maps
$$\gamma_{n,t} : H^{1}(\FS/\Fn,\At) \to \prod_{w \mid v \in \Sigma}
\Hs(\Fnw,\At)$$
have finite cokernel bounded independent of $n$.
Indeed, taking the limit in $n$ we then find that
$$\gamma_{t} : H^{1}(\FS/\Fi,\At) \to \prod_{w \mid v \in \Sigma}
\Hs(\Fiw,\At)$$
has finite cokernel; 
here $\Hs(\Fiw,\At)$ is defined analogously to $\Hs(\Fnw,\At)$ above.  Since
$\At \cong A$ over $\Fi$, one sees immediately that the cokernel of $\gamma$
is a quotient of the cokernel of $\gamma_{t}$, so that this does indeed
suffice.

We now select an appropriate value of $t$.  We claim that for all but
finitely many $t \in \Zp$ we have:
\renewcommand{\theenumi}{\roman{enumi}}
\begin{enumerate}
\item $\ker \gamma_{n,t}$ is finite for all $n$;
\item $H^{0}(\Fnw,\At)$ and $H^{0}(\Fnw,\At^{*})$ are finite for all 
$w \mid v \in \Sigma$, $v \nmid p\infty$, and all $n$;
\item $H^{0}(\Fnw,\At/\Atw^{\crit})$ and $H^{0}(\Fnw,(\At/\Atw^{\crit})^{*})$
are finite for all $w \mid p$ and all $n$;
\item $H^{0}(\Fnw,(\Atw^{\crit})^{*})$ is finite for all $w \mid p$ and all $n$;
\end{enumerate}
This is clear for the latter three conditions as $\kappa|_{\Gw}$ 
has infinite order for each of the finitely many non-archimedean
places of $\Fi$ lying
over places in $\Sigma$.
For the first condition, consider the restriction map
$$\ker \gamma_{n,t} \to \Selc(\Fi,\At)^{\Gamma^{p^{n}}}$$
where $\Gamma$ denotes $\Gal(\Fi/F)$.  The kernel of this map lies
in $H^{1}(\Fi/\Fn,\At^{\GFi})$, which has the same corank as
$H^{0}(\Fn,\At)$ and thus is finite for almost all $t$.
Since $\Selc(\Fi,A)$ is $\Lambda$-cotorsion by assumption,
$\Selc(\Fi,\At)^{\Gamma^{p^{n}}}$ is finite for all $n$ for almost all $t$;
any $t$ satisfying the latter two conditions satisfies (i).

Fix such a $t$.
By (ii) and local duality we have that
$H^{2}(\Fnw,\At)$ is finite for all $w$ dividing $v
\in \Sigma$, $w \nmid p$, and all $n$.
By (ii) and the local Euler characteristic formula, 
it follows that $H^{1}(\Fnw,\At)$ is finite for all such $w$ and all $n$.
By an analogous argument using (iii) one sees that
$H^{1}(\Fnw,\At/\At^{\crit})$ has corank $[\Fnw:\Qp]c(V)
= [\Fv:\Qp]c(V)p^{n}$ 
for all $w$ dividing $v$ dividing $p$.  By (iv) and another
application of local duality this implies that
$\Hs(\Fnw,\At)$ has corank $[\Fv:\Qp]c(V)p^{n}$ as well.

The Poitou--Tate global duality sequence yields an exact sequence
$$0 \to \ker \gamma_{n,t} \to H^{1}(\FS/\Fn,\At) \overset{\gamma_{n,t}}{\too}
\prod_{w \mid v \in \Sigma} \Hs(\Fnw,\At) \to \\
H_{1,n} \to H_{2,n} \to 0$$
where $H_{1,n}$ is dual to a subgroup of $H^{1}(\FS/\Fn,\At^{*})$ and
$H_{2,n}$ is a subgroup of $H^{2}(\FS/\Fn,\At)$.
By the global Euler characteristic formula and the definition of $c(V)$
we have
$$\corank_{\O} H^{1}(\FS/\Fn,\At) = [F:\Q]c(V)p^{n} + 
\corank_{\O} H^{2}(\FS/\Fn,\At).$$
By our local computations above we see that the target of $\gamma_{n,t}$
has corank $[F:\Q]c(V)p^{n}$.  
Since $\ker \gamma_{n,t}$ is finite by (i), it follows
that $H^{1}(\FS/\Fn,\At)$ has corank $[F:\Q]c(V)p^{n}$ and that
$H^{2}(\FS/\Fn,\At)$ is finite.  
In fact, $H^{2}(\FS/\Fi,\At)$ must therefore vanish since by 
\cite[Proposition 4]{Greenberg1} it is cofree over $\Lambda$.
Thus $H_{2,n}$ vanishes, so that $H_{1,n} = \coker \gamma_{n,t}$ 
is finite as well.

In particular, $\coker \gamma_{n,t}$ must be dual to a subgroup
of $H^{1}(\FS/\Fn,\At^{*})_{\tors}$.  This latter group is simply the
kernel of
$$H^{1}(\FS/\Fn,\At^{*}) \to H^{1}(\FS/\Fn,\At^{*} \otimes K)$$
and thus equals the image of $H^{0}(\Fn,\At^{*} \otimes K/\O)$.
Thus for any $n$ we can identify the dual of $\coker \gamma_{n,t}$ with a
subgroup of $H^{0}(\Fi,\At^{*} \otimes K/\O) = H^{0}(\Fi, A^{*} \otimes K/\O)$.
This latter
group is finite by assumption, so that
$\coker \gamma_{n,t}$ is bounded independent of
$n$; the first part of the proposition follows.

We turn now to (2).  Let $t$ be as above, subject to the additional
hypothesis that
$$H^{1}(\FS/\Fi,\At)/H^{1}(\FS/\Fi,\At)_{\Lambda-\text{div}}$$
has finite $\Gamma$-covariants.  (This is certainly possible since
$$H^{1}(\FS/\Fi,A)/H^{1}(\FS/\Fi,A)_{\Lambda-\text{div}}$$
is visibly $\Lambda$-cotorsion.)
Since the cokernel of the injection
$$\ker \gamma_{t} \inj \Selc(\Fi,A)$$
is divisible, 
to show that $\Selc(\Fi,A)$ has no proper $\Lambda$-submodules of finite
index it suffices to show the same for $\ker \gamma_{t}$.
In fact, by
the structure theory of $\Lambda$-modules for this it suffices to show that
$(\ker \gamma_{t})_{\Gamma} = 0$.

Since $\Gamma$ is pro-$p$, 
the assumption that $H^{0}(F,A^{*} \otimes K/\O)=0$ implies that
$$H^{0}(\Fi,A^{*} \otimes K/\O)= H^{0}(\Fi,\At^{*} \otimes K/\O) = 0.$$
In particular, the group $H_{1,0}$
above in fact vanishes, so that the map
$$H^{1}(\FS/F,\At) \to \prod_{w \mid v \in \Sigma} \Hs(\Fw,\At)$$
is surjective.  It follows from this and the fact that $\Gamma$
has cohomological dimension one that the map
$$H^{1}(\FS/\Fi,\At)^{\Gamma} \to \prod_{w \mid 
v \in \Sigma} \Hs(\Fiw,\At)^{\Gamma}$$
is surjective as well.  This in turn implies that
the map
$$(\ker \gamma_{t})_{\Gamma} \to H^{1}(\FS/\Fi,\At)_{\Gamma}$$
is injective.
The latter group is finite by our last assumption on $t$.
However, $H^{1}(\FS/\Fi,\At)$ has no proper $\Lambda$-submodules 
of finite index
by \cite[Proposition 5]{Greenberg1}.
Thus $H^{1}(\FS/\Fi,\At)_{\Gamma}$ must in fact vanish, so that
$(\ker \gamma_{t})_{\Gamma}$ vanishes as well, as desired.
\end{proof}

\subsection{Iwasawa invariants}

Fix a uniformizer $\pi$ of $\O$ and set $k=\O/\pi$.
Assume now that $A$ is a torsion quotient of a
critical nearly ordinary Galois representation $V$; that
$\Selc(\Fi,A)$ is $\Lambda$-cotorsion; and that
$H^{0}(F,A^{*} \otimes K/\O)=0$.
We define $\mu(A)$ to be the least $n \geq 0$ such that
$$\Selc(\Fi,A)[\pi^{n+1}]/\Selc(\Fi,A)[\pi^{n}]$$ 
is finite-dimensional over $k$.  
We define $\lambda(A)$ as the $\O$-corank of $\Selc(\Fi,A)$.
Proposition~\ref{prop:loc} yields the following result.

\begin{corollary} \label{cor:ii}
Let $A$ be a torsion quotient of a critical nearly ordinary Galois
representation $V$ such that $\Selc(\Fi,A)$ is $\Lambda$-cotorsion and
$H^{0}(F,A^{*} \otimes K/\O)=0$.  Then $\mu(A) = 0$ if
and only if
$\Selc(\Fi,A)[\pi]$ is finite-dimensional.   If $\mu(A)=0$, then
$\Selc(\Fi,A)$ is a divisible $\O$-module and
$$\lambda(A) = \dim_{k} \Selc(\Fi,A)[\pi].$$
\end {corollary}
\begin{proof}
The first statement is immediate from the definition.  If this is
the case, then it follows from (2) of Proposition~\ref{prop:loc} that
$\Selc(\Fi,A)$ is divisible: indeed, if $\mu(A) = 0$ then
the maximal divisible subgroup
of $\Selc(\Fi,A)$ has finite index, so that by the proposition it must
coincide with $\Selc(\Fi,A)$.  By definition we thus have
$$\Selc(\Fi,A) \cong (K/\O)^{\lambda(A)};$$
the corollary follows.
\end{proof}

\section{Residual Selmer groups} \label{s2}

Corollary~\ref{cor:ii} suggests that it should be possible to study
the Iwasawa invariants of a critical torsion quotient $A$ by studying
a corresponding residual Selmer group.  However, this Selmer group must
depend not
only on the $\pi$-torsion of $A$ but also on the ramification of $A$.
In this section we model the different possible Selmer groups using a version
of Mazur's notion of a finite/singular structure.

\subsection{Ordinary residual Galois representations} \label{s21}

Let $F$ be a number field as before and
let $k$ be a finite field of odd characteristic $p$.  An {\it ordinary
Galois representation} over $k$ consists of a finite-dimensional
$k$-vector space $\Ab$ endowed with a $k$-linear action of $\GF$ and
choices of $\Gv$-stable complete flags
$$0 = \Ab_{v}^{0} \subsetneq \Ab_{v}^{1} \subsetneq \cdots
\subsetneq \Ab_{v}^{n} = \Ab$$
for each place $v$ dividing $p$.

Let $\Ab$ be an ordinary Galois representation of dimension $n$ over $k$.
For a real place $v$ of $F$ let $\dvm(\Ab)$
denote the $K$-dimension of the subspace of $V$ on which complex
conjugation at $v$ acts by
$-1$.
We say that $\Ab$ is {\it critical} if
$$\sum_{v \text{~real}} \dvm(\Ab) + \sum_{v \text{~complex}} n$$
is divisible by $[F:\Q]$, in which case we write $c(\Ab)$ for the quotient
and set $\Ab_{v}^{\crit} = \Ab_{v}^{n-c(\Ab)}$ for each place
$v$ dividing $p$.

\subsection{Finite/singular structures} \label{s22}

Let $\Ab$ be a critical ordinary Galois representation over
$k$.

\begin{defn}
A {\it finite/singular} structure $\S$ on $\Ab$ is a choice of
$k$-subspaces
$$\HfS(\Fiw,\Ab) \subseteq H^{1}(\Gw,\Ab)$$
for every place $w$ of $\Fi$ subject to the restrictions:
\begin{enumerate}
\item $\HfS(\Fiw,\Ab) = \ker H^{1}(\Gw,\Ab) \to H^{1}(\Iw,\Ab_{w}^{\crit})$
for any place $w$ dividing $p$;
\item $\HfS(\Fiw,\Ab) = 0$ for almost all $w$;
\item $\HfS(\Fiw,\Ab)$ and $\HfS(F_{\infty,w'},\Ab)$ coincide under
the canonical isomorphism $H^{1}(\Gw,\Ab) \cong H^{1}(G_{w'},\Ab)$
for any places $w$ and $w'$ dividing the same place $v \nmid p$ of $F$.
\end{enumerate}
(Note that $H^{1}(\Gw,\Ab)=0$ for all archimedean places $w$, so that we
may safely ignore archimedean places below.)
\end{defn}

Fix a finite/singular structure $\S$ on $\Delta$.
We say that a finite set $\Sigma$ of places of $F$ is {\it sufficiently
large} for $\S$ if it contains all places at which $\Ab$ is ramified, all
archimedean places, all places dividing $p$, and all places $v$ such that
$\HfS(\Fiw,\Ab)$ is non-zero for some $w$ dividing $v$.
For a place $v$ of $F$ we set
$$\delta_{\S,v}(\Ab) = \sum_{w \mid v}
\dim_{k} \HfS(\Fiw,\Ab)$$
and for a place $w$ of $\Fi$ we set
$$\HsS(\Fiw,\Ab) = H^{1}(\Gw,\Ab)/\HfS(\Fiw,\Ab).$$
(Note that $H^{1}(\Gw,\Ab)$ is always finite dimensional over $k$,
so that $\delta_{\S,v}(\Ab)$ is finite.)
We define the {\it $\S$-Selmer group} of $\Delta$ by
\begin{align}
\SelS(\Fi,\Ab) &= \ker H^{1}(\Fi,\Ab) \to \prod_{w} \HsS(\Fiw,\Ab) \notag \\
&= \ker H^{1}(\FS/\Fi,\Ab) \to \prod_{w \mid v \in \Sigma} \HsS(\Fiw,\Ab)
\end{align}
for any finite set of places $\Sigma$ of $F$ 
which is sufficiently large for $\S$.

We say that $\mu(\Delta)=0$ if $H^{1}(\FS/\Fi,\Ab)$ is finite dimensional;
this is independent of the choice of sufficiently large $\Sigma$.
If this is the case, then $\SelS(\Fi,\Ab)$ is finite dimensional for
any $\S$, and we set
$$\lambda_{\S}(\Ab) = \dim_{k} \SelS(\Fi,\Ab).$$

\begin{example}
The {\it minimal structure} $\Sm$ is given by
$$\HfSm(\Fiw,\Ab) = 0$$
for $w$ not dividing $p$.  (The condition at $w$ dividing $p$ is
fixed by definition.)
The corresponding {\it minimal Selmer group}
$\SelSm(\Fi,\Ab)$ is contained in every other Selmer group of $\Ab$.
\end{example}

\begin{example}
Let $K$ be a finite extension of $\Qp$ with residue field $k$; we write
$\O$ for the ring of integers and $\pi$ for a fixed choice of uniformizer.
Let $A$ be a torsion quotient of a
critical nearly ordinary Galois representation $V$ over $K$.
The $\pi$-torsion $A[\pi]$ inherits an obvious structure of 
critical ordinary
Galois representation over $k$.
We define a structure $\S(A)$ on $A[\pi]$ by setting
$$\HfSA(\Fiw,A[\pi]) = \ker H^{1}(\Gw,A[\pi]) \to \Hsc(\Fiw,A)$$
for every place $w$ of $\Fi$.
Note that for $w$ not dividing $p$ we have
$$\HfSA(\Fiw,A[\pi]) = \im A^{\Gw}/\pi \inj H^{1}(\Gw,A[\pi]).$$
In order to prove that $\S(A)$ is a finite/singular structure on $A[\pi]$
we need a simple yet crucial lemma on local Galois invariants.

\begin{lemma} \label{lemma:div}
Let $A$ be a torsion quotient of a Galois representation $V$ over $K$.
Let $v$ be a place of $F$ not dividing $p$ such that
$V$ and $A[\pi]$ have the same
Artin conductor at $v$.  Then $A^{\Gw}$ is $\O$-divisible for every place
$w$ of $\Fi$ dividing $v$.
\end{lemma}
\begin{proof}
Fix a place $w$ of $\Fi$ dividing $v$.
As the Swan conductor is invariant under reduction
(see \cite{Livne} for example), it follows
from the hypothesis and the definition of the Artin conductor that
$$\dim_{K} V^{\Iw} = \dim_{k} A[\pi]^{\Iw}.$$
(Note that $\Iw = \Iv$ since $\Fi/F$ is unramified at $v$.)
Since $A^{\Iw}$ has $\O$-corank equal to $\dim_{K} V^{\Iw}$, we
conclude that $A^{\Iw}$ is $\O$-divisible.
The $\Gw/\Iw$-invariants of an $\O$-divisible module are still 
$\O$-divisible, so that the lemma follows from this.
\end{proof}

\begin{proposition} \label{prop:structure}
Let $A$ be a torsion quotient of a Galois representation $V$ over $K$.
If $H^{0}(\Iw,A/A_{w}^{\crit})$ is $\O$-divisible
for each place $w$ dividing $p$, then
the structure $\S(A)$ is a finite/singular structure on $A[\pi]$.
\end{proposition}

The divisibility hypothesis above appears to be essential to our method.
In the case that $V$ arises from a modular form $f$, it corresponds
to the assumption that $f$ is a twist of an ordinary modular form (in the
usual sense that the $p^{\text{th}}$ Fourier coefficient is prime to $\pi$)
by a power of the Teichm\"uller character.

\begin{proof}
It follows from Lemma~\ref{lemma:div} that $\HfSA(\Fiw,A[\pi])=0$ for
almost all $w$.  The coincidence of $\HfSA(\Fiw,A[\pi])$ and
$\HfSA(F_{\infty,w'},A[\pi])$ for $w$ and $w'$ dividing $v \nmid p$
is immediate from the definition.
Finally, the verification of the conditions at places $w$
dividing $p$
is a simple diagram chase using the fact that
the divisibility of $H^{0}(\Iw,A/A_{w}^{\crit})$ implies that
$$H^{1}(\Iw,A/A_{w}^{\crit}[\pi]) \to H^{1}(\Iw,A/A_{w}^{\crit})$$
is injective
\end{proof}

\begin{corollary} \label{cor:res}
Let $A$ be a torsion quotient of a Galois representation $V$ over $K$
such that $H^{0}(\Iw,A/A_{w}^{\crit})$ is $\O$-divisible
for each place $w$ dividing $p$.
If $H^{0}(\Fi,A)$ is $\O$-divisible, then the natural map
$$\SelSA(\Fi,A[\pi]) \to \Selc(\Fi,A)[\pi]$$
is an isomorphism.  In particular, $\mu(A)=0$ if and only if
$\mu(A[\pi])=0$.  If this is the case, then
$\lambda(A) = \lambda_{\S(A)}(A[\pi])$.
\end{corollary}
\begin{proof}
The identification of Selmer groups is immediate from the
injectivity of
$$H^{1}(\Fi,A[\pi]) \to H^{1}(\Fi,A)$$
and the definition of $\S(A)$.  The relation between the Iwasawa invariants
now follows from the definitions and Corollary~\ref{cor:ii}.
\end{proof}

In particular, the corollary implies that
knowledge of the residual representation $A[\pi]$
and finite/singular structure $\S(A)$ (which depends only on
$A[\pi]$ and the ramification of $A$) determines the Iwasawa invariants
of $\Selc(\Fi,A)$.  In the next section we will use
Proposition~\ref{prop:loc}
to give a more precise description of this relation.

\end{example}

We remark that one can consider induced structures as above for lifts of
$\Ab$ to representations over
more general complete local noetherian rings with residue field
$k$.  In particular, one can then compare the structures induced from the
ring and its quotients.  For example, one can show that induced structures
are constant on families with constant ramification in an appropriate sense.
(In the case of Hida families, by this we mean a branch of the Hida family
with all crossing points removed.)  We will not pursue this point of view
any further here.

\subsection{Variation of structure} \label{s23}

Fix a critical ordinary Galois representation $\Ab$
over $k$.  Let $K$ be a finite extension of $\Qp$ with residue field $k$,
ring of integers $\O$ and uniformizer $\pi$.
A {\it lift} of $\Ab$ over $K$ is a
torsion quotient $A$ of a 
nearly ordinary Galois representation $V$ over $K$ such
that $H^{0}(\Iw,A/A_{w}^{\crit})$ is $\O$-divisible for all $w$ dividing $p$,
together with  an isomorphism $\Ab \cong A[\pi]$ of nearly ordinary
Galois representations.  Note that $V$ is then necessarily 
critical since $p$ is odd and thus
$\dvm(V) = \dvm(\Ab)$ for all real places $v$.  
We say that a finite/singular structure $\S$ on $\Delta$ is
{\it induced} if there is a lift $A$ of $\Ab$ such
that the isomorphism $\Ab \cong A[\pi]$ identifies $\S$ with $\S(A)$.
We say that $\S$ is {\it properly induced} if such an $A$ can be
chosen so that $\Selc(\Fi,A)$ is $\Lambda$-cotorsion.

The key result for the analysis of residual Selmer groups is the following
proposition.

\begin{proposition} \label{prop:exact}
Assume that $\mu(A) = 0$ and $H^{0}(F,\Ab)=H^{0}(F,\Ab^{*})=0$.
If $\S$ is properly induced, then the sequence
$$0 \to \SelS(\Fi,\Ab) \to H^{1}(\FS/\Fi,\Ab) \to
\prod_{w \mid v \in \Sigma} \HsS(\Fiw,\Ab) \to 0$$
is exact for any finite set of places $\Sigma$ of $F$
sufficiently large for $\S$.
\end{proposition}
\begin{proof}
Fix a proper lift $A$ of $\Ab$ which
identifies $\S(A)$ with $\S$.
Note that the hypotheses imply that
$$H^{0}(F,A) = H^{0}(F,A^{*} \otimes K/\O) = 0.$$
Consider the exact sequence
$$0 \to \Selc(\Fi,A) \to H^{1}(\FS/\Fi,A) \to \prod_{w \mid v \in \Sigma}
\Hsc(\Fiw,A) \to 0$$
of Proposition~\ref{prop:loc}.  Since $\Selc(\Fi,A)$
is $\O$-divisible by Corollary~\ref{cor:ii}, the $\pi$-torsion of
this sequence is an exact sequence
$$0 \to \Selc(\Fi,A)[\pi] \to H^{1}(\FS/\Fi,A)[\pi] \to \prod_{w 
\mid v \in \Sigma}
\Hsc(\Fiw,A)[\pi] \to 0.$$
Since $H^{0}(\Fi,A)=0$ is divisible,
it follows easily from Corollary~\ref{cor:res} and the definitions
that this sequence identifies with the desired sequence.
\end{proof}

The above proof rests entirely on the crutch of a proper lift of $\Ab$.
We do not know how to approach this problem purely via the residual
representation $\Ab$.

\begin{corollary} \label{cor:quant}
Assume that $\mu(\Ab)=0$ and $H^{0}(F,\Ab)=
H^{0}(F,\Delta^{*})=0$.  The quantity
$$\lambda_{\S}(\Ab) - \sum_{v \nmid p} \delta_{\S,v}(\Ab)$$
is independent of the choice of properly induced structure $\S$.
\end{corollary}
\begin{proof}
Let $\S_{1}$ and $\S_{2}$ be properly induced structures on $\Ab$ and
choose $\Sigma$ which is sufficiently large for both.  Then
it follows from Proposition~\ref{prop:exact} and the agreement
of the local conditions at places dividing $p$ that
$$\lambda_{\S_{1}}(\Ab) + \sum_{\substack{w \mid v \in \Sigma\\v \nmid p}}
\dim H^{1}_{s,\S_{1}}(\Fiw,\Ab) =
\lambda_{\S_{2}}(\Ab) + \sum_{\substack{w \mid v \in \Sigma\\
v \nmid p}} \dim H^{1}_{s,\S_{2}}(\Fiw,\Ab).$$
Since
$$\delta_{\S_{i},v}(\Ab) = \sum_{w \mid v} \dim H^{1}(\Fiw,\Ab) - \dim
H^{1}_{s,\S_{i}}(\Fiw,\Ab)$$
the corollary follows from this.
\end{proof}

Note that when the minimal structure $\Sm$ itself is properly induced the
above difference is simply $\lambda_{\Sm}(\Ab)$.
The question of whether or not the minimal structure
$\Sm$ is induced is intimately related to level lowering.
More precisely, we have the following result, which follows
immediately from Lemma~\ref{lemma:div}.

\begin{proposition} \label{lemma:min}
Let $A$ be a lift of $\Ab$ over $K$.
If the Artin conductor of $\Ab$ equals
the Artin conductor of $V$, then the isomorphism $\Ab \cong A[\pi]$
identifies $\Sm$ with $\S(A)$.
\end{proposition}

\subsection{Existence of structures} \label{s24}

Proposition~\ref{prop:exact} gives very precise control over properly
induced structures.  It thus becomes an interesting question to
determine which structures actually occur in this way; that is,
for which structures $\S$ on $\Ab$ do there exist torsion quotients
$A$ of critical nearly ordinary
Galois representations $V$ with $\Lambda$-cotorsion Selmer group and an
isomorphism $\Ab \cong A[\pi]$ identifying $\S$ with $\S(A)$?

The local part of this question is not difficult: given $\Ab$, a place $w$ of
$\Fi$ dividing a place $v$ of $F$, and a subspace
$$H \subseteq H^{1}(\Gw,\Ab)$$
it is a straightforward calculation
to determine if there is a $\Gv$-representation $A$
lifting $\Ab$ such that
$$ H = \im A^{\Gw}/\pi \to H^{1}(\Gw,\Ab).$$
Much more difficult is the amalgamation of this local information into
a global structure.  This latter question is intimately connected with
level raising in the sense of \cite{DiamondTaylor2}.  To illustrate this
connection we give the following result for two-dimensional 
modular representations of tame level one; one can prove similar results
for higher levels, but we focus on this case for simplicity.

Let $\Ab$ be an absolutely irreducible ordinary 
two-dimensional modular Galois representation over $\Q$ unramified
away from $p$.
Assume also that the $G_{p}$-representation
$\Ab/\Ab_{p}^{1}$ is unramified.
Let $\S$ be a finite/singular structure on
$\Ab$.  Consider the two conditions:
\begin{enumerate}
\item $\HfS(\Qiw,\Ab)=0$ for any place $w \nmid p$ dividing a prime
$\ell \not\equiv 1 \pmod{p}$ such that
$$\Ab|_{\Gl} \cong \left( \begin{array}{cc}
\omega & * \\ 0 & \chi \end{array} \right)$$
with $\omega$ the Teichm\"uller character and $\chi \neq 1$.
(Note that $\Ab|_{\Gl}$ is unramified by assumption and thus
will be split in this case unless $\chi = \omega$.)
\item $\dim \HfS(\Qiw,\Ab) \neq 1$ for any place $w$ dividing a prime
$\ell \equiv 1 \pmod{p}$ such that
$\Ab|_{\Gl}$ is trivial.
\end{enumerate}

\begin{proposition} \label{prop:levelraising}
Let $\S$ be a finite/singular structure on $\Delta$.
If $\S$ is induced, then (1) holds.  If (1) and (2) hold,
then $\S$ is properly induced.
\end{proposition}

It is an interesting question as to whether or not condition
(1) alone suffices to determine if $\S$ is induced.  This reduces to
the question of whether or not for a prime $\ell \equiv 1 \pmod{p}$ 
such that $\Ab|_{\Gl}$ is trivial one can find modular lifts of
$\Ab$ which are special (of level $\ell$) at $\ell$ and with the
unramified line lifting an arbitrary choice of line in $\Ab$.
It is not clear to
this author if the techniques of \cite{DiamondTaylor2} can be
modified to answer this question.

In the proof we will frequently use the fact that
$H^{1}(\Gw,\Ab)$ (resp.\ $H^{1}(\Gw,A)[\pi]$) has
$k$-dimension equal to the multiplicity of $\omega$ in
the inertia coinvariants $\Ab_{\Iw}$ (resp.\ $V_{\Iw}$);
see \cite[Section 2]{GreenbergVatsal}.

\begin{proof}
Suppose first that $\S$ is induced and let $w$ be a place dividing a
prime $\ell$ with
$\Ab|_{\Gl}$  as in (1).  Let $A$ be a torsion quotient of an
ordinary Galois representation $V$ lifting $\S$ and
consider the corresponding exact sequence
$$0 \to A^{\Gw}/\pi \to H^{1}(\Gw,\Ab) \to H^{1}(\Gw,A)[\pi] \to 0.$$
If $A$ is unramified at $\ell$,
then $A^{\Gw}$ is divisible, so that we must have $\HfS(\Qiw,\Ab)=0$.
If $A$ is ramified at $\ell$, then since $\Ab$ is unramified at $\ell$
yet $\ell \not\equiv 1 \pmod{p}$ and $\chi \neq 1$, we must have
$\chi = \omega^{2}$.  
(See \cite[Section 1]{DiamondTaylor2}.)
But then $\ell \not\equiv -1 \pmod{p}$ (for that
would force $\chi = 1$) so that $V$ must be special at $\ell$:
that is, the Galois action on $V$ is given by
$$\left( \begin{array}{cc} \omega^{2} & * \\ 0 & \omega \end{array} \right)$$
with $*$ ramified but trivial modulo $\pi$.  Then both
$H^{1}(\Gw,\Ab)$ and $H^{1}(\Gw,A)[\pi]$ are one-dimensional,
so that we still have $\HfS(\Qiw,\Ab)=0$.

Next suppose that (1) and (2) are satisfied and fix a prime $\ell$
different from $p$.  We claim that there exists a
two-dimensional $K$-vector space $V_{\ell}$ with
a $K$-linear action of $\Gl$ such that
$$\sum_{w \mid \ell} \dim_{k} H^{1}(\Gw,\Ab) - \dim_{K} H^{1}(\Gw,V_{\ell}) =
\delta_{\S,\ell}(\Ab).$$
Indeed, this is straightforward from the discussion of
\cite[Section 1]{DiamondTaylor2}.  Specifically, we may take $V_{\ell}$
unramified if $\delta_{\S,\ell}(\Ab)=0$.
If $\delta_{\S,\ell}(\Ab) \neq 0$ and
$\ell \not\equiv 1 \pmod{p}$, then by (1) we must have
$$\Ab|_{\Gl} \cong
\left( \begin{array}{cc} \omega & 0 \\ 0 & 1 \end{array} \right).$$
In this case one simply takes $V_{\ell}$ to be the special representation
$$\left( \begin{array}{cc} \cyc & * \\ 0 & 1 \end{array} \right)$$
with $*$ trivial modulo $\pi$.
If $\delta_{\S,\ell}(\Ab) \neq 0$ and $\ell \equiv 1 \pmod{p}$,
then one may take $V_{\ell}$ to be a ramified principal series lifting
of $\Ab$.

After possibly enlarging $K$,
we may now apply \cite[Theorem 1]{DiamondTaylor2} to obtain a
torsion quotient $A$ of a (necessarily ordinary) modular representation
$V$ together with an isomorphism $\Ab \cong A[\pi]$ and
isomorphisms $V|_{\Il} \cong V_{\ell}$ for every $\ell \neq p$.
Since $\dim_{K} H^{1}(\Gw,V)$ depends only on $\Ab$ and $V|_{\Il}$,
we have
$$\sum_{w \mid \ell} \dim_{k} H^{1}(\Gw,\Ab) - \dim_{K} H^{1}(\Gw,V) =
\delta_{\S,\ell}(\Ab)$$
for all places $w$ dividing $\ell \neq p$.
Under the assumption (2) this is in fact enough to determine
$\HfS(\Qiw,\Ab)$, so that we must have $\S(A)=\S$.
Finally, since $V$ is modular it is known by work of Kato
that $\Selc(\Qi,A)$ is $\Lambda$-cotorsion, so that $\S$ is
properly induced, as desired.
\end{proof}

\section{Galois deformations} \label{s3}

\subsection{Deformations} \label{s31}

Fix a finite extension $K_{0}$ of $\Qp$ with ring of integers $\O_{0}$
and residue field $k$.
Let $\G$ denote a split reductive algebraic group over $\O_{0}$ such that
the center $\ZG$ of $\G$ is smooth over $\O_{0}$.  
For each place $v$ dividing $p$ fix a Borel subgroup $\B_{v}$ of $\G$.
Let
$$\rhob : \GF \to \G(k)$$
be a continuous residual Galois representation taking values in
$\G$.  We assume:
\begin{enumerate}
\item the identity component of the centralizer of $\rhob$ is equal to
the identity component of  $\ZG \times k$;
\item $\rhob(\Gv)$ lies in $\B_{v}(k)$ for every place $v$ dividing $p$;
\item $\rhob$ is {\it regular}: for any place $v$ dividing $p$ the
invariants $H^{0}(\Gv,\Glie/\Blie_{v})$ vanish.  (Here $\Glie$ and
$\Blie_{v}$ are the Lie algebras of $\G$ and $\B_{v}$ endowed with adjoint
$\Gv$-action.)
\end{enumerate}

A {\it nearly ordinary lift} $\rho$ 
of $\rhob$ to a local $\O_{0}$-algebra $\O$ with
residue field $k$ is a continuous representation
$$\rho : \GF \to \G(\O),$$
ramified at finitely many places, such that the composition
$$\GF \overset{\rho}{\to} \G(\O) \to \G(k)$$
is equal to $\rhob$ and such that for each place $v$ dividing $p$ there is
$$g_{v} \in \Ghat(\O) := \ker \G(\O) \to \G(k)$$
such that
$$g_{v} \cdot \rho(\Gv) \cdot g_{v}^{-1} \subseteq \B_{v}(\O).$$
We consider two such liftings equivalent if one can be conjugated to
the other via some element of $\Ghat(\O)$;
a {\it nearly ordinary 
deformation} of $\rhob$ is an equivalence class of liftings.

For a local $\O_{0}$-algebra $\O$ with residue field $k$
we write $\H(\O)$ for the set of nearly ordinary
deformations of $\rhob$ to $\O$.  If $\Sigma$ is a finite set of places
sufficiently large for $\rhob$, we write
$\H_{\Sigma}(\O)$ for the subset of $\H(\O)$ of deformations unramified
away from $\Sigma$.  By \cite[Proposition 6.2]{Tilouine}
the functor $\H_{\Sigma}$ is representable: there is a local $\O_{0}$-algebra
$R_{\Sigma}^{\no}$ with residue field $k$ such that there is a bijection
between $\H(\O)$ and $\Hom(R_{\Sigma}^{\no},\O)$ for any $\O$ as above.
(We note that $\Hom(R_{\Sigma}^{\no},\O)$ is computed in the category of
inverse limits of artinian local $\O_{0}$-algebras with residue field $k$;
such morphisms are required to be local homomorphisms
inducing the identity map on $k$.)

\subsection{Selmer groups} \label{s32}

Fix an algebraic representation
$$r : \G \to \GL_{n}$$
for some $n$.  For each place $v$ dividing $p$ fix also a Borel
subgroup $\tilde{\B}_{v}$ of $\GL_{n}$ containing $r(\B_{v})$.
Let $\Ab(\rhob,r)$ denote the representation space for
$$\GF \overset{\rhob}{\too} \G(k) \overset{r}{\too} \GL_{n}(k).$$
We endow $\Ab(\rhob,r)$ with the structure of 
ordinary Galois representation by letting
$$0 = \Ab_{v}^{0} \subsetneq \Ab_{v}^{1} \subsetneq \cdots
\subsetneq \Ab_{v}^{n} = \Ab(\rho,r)$$
be the complete flag associated to the Borel subgroup $\tilde{\B}_{v}$ for each
$v$ dividing $p$; since $r \circ \rhob(\Gv) \subseteq \tilde{\B}_{v}(k)$,
this flag is indeed $\Gv$-stable.  We say that the pair $(\rhob,r)$ is
{\it critical} if $\Ab(\rhob,r)$ is critical.

Let $\O$ be the ring of integers of a finite totally ramified extension $K$
of $K_{0}$ and fix $\rho \in \H(\O)$.
Define a Galois representation
$\Tr$ as the representation space of the composition
$$\GF \overset{\rho}{\to} \G(\O) \overset{r}{\to}
\GL_{n}(\O).$$
The isomorphism class of $\Tr$ depends only on the deformation
class of $\rho$, so that we may speak of the Galois representation
associated to $\rho$.  Set $\Vr = \Tr \otimes_{\O} K$ and
$\Ar = \Vr/\Tr$.

We claim that we may endow $\Vr$ with a canonical structure of
nearly ordinary Galois representation.  Indeed, let $v$ be a place dividing
$p$ and fix $g_{v} \in \Ghat(\O)$ such that
$$g_{v} \cdot \rho(\Gv) \cdot g_{v}^{-1} \subseteq \B_{v}(\O).$$
Define a complete flag
$$0 = \Vrv^{0} \subsetneq \Vrv^{1} \subsetneq \cdots \subsetneq
\Vrv^{n} = \Vr$$
as the $g_{v}^{-1}$ conjugate of the complete flag associated to
the Borel subgroup $\tilde{\B}_{v}$ of $\GL_{n}$.  It is clear that
this flag is $\Gv$-stable.  It is also independent of the choice
of $g_{v}$ as above: by \cite[Claim of p.\ 49]{Tilouine} (using
the regularity assumption) any other
$g_{v}' \in \Ghat(\O)$ such that 
$$g_{v}' \cdot \rho(\Gv) \cdot
g_{v}'{}^{-1} \subseteq \B_{v}(\O)$$
must differ from $g_{v}$ by an element of
$$\widehat{\B}_{v}(\O) := \ker \, \B_{v}(\O) \to \B_{v}(k)$$
and thus yields the same complete flag.  
We always regard $\Vr$ as a nearly ordinary
Galois representation via these choices of complete flags.

Assume that $(\rhob,r)$ is critical and
define the Selmer group $\Sel(\Fi,\rho,r)$ by
$$\Sel(\Fi,\rho,r) := \Selc(\Fi,\Ar).$$
(Of course, this Selmer group also depends on the choices of Borel
subgroups $\tilde{\B}_{v}$; we omit these choices from the
notation.)  When $\Sel(\Fi,\rho,r)$ is $\Lambda_{\O}$-cotorsion we
let $\mu(\rho,r)$ and
$\lambda(\rho,r)$ denote the $\mu$ and $\lambda$-invariants of
$\Sel(\Fi,\rho,r)$.

Let $\rho_{1} \in \H(\O_{1})$ and $\rho_{2} \in \H(\O_{2})$ be
two deformations as above.  
We say that $\rho_{1} \in \H(\O_{1})$ and
$\rho_{2} \in \H(\O_{2})$ are {\it stably equivalent} if there is
a finite totally ramified extension $K_{3}$ of $K_{0}$, with ring of integers
$\O_{3}$, and injective morphisms
$\sigma_{1} : \O_{1} \to \O_{3}$ and 
$\sigma_{2} : \O_{2} \to \O_{3}$
such that $\sigma_{1} \circ \rho_{1}$ and $\sigma_{2} \circ \rho_{2}$
are equivalent deformations over $\O_{3}$.
We simply write $\H$ for the set of stable equivalence classes of
nearly ordinary deformations of $\rhob$; since $\GF$ is compact it can
also be interpreted as the set of nearly ordinary deformations of $\rhob$
to the ring of integers of the maximal totally ramified extension of $K_{0}$.
Note that as defined above the vanishing of $\mu(\rho,r)$ and the value of
$\lambda(\rho,r)$ depend only
on the stable equivalence class of $\rho$

Let $\H(r)$ denote the subset of $\rho \in \H$ such that
for some (or equivalently for all) choice of representative $\rho \in \H(\O)$,
$\Sel(\Fi,\rho,r)$ is $\Lambda_{\O}$-cotorsion and such that
$H^{0}(\Fiw,A_{\rho}/A_{\rho,w}^{\crit})$ is $\O$-divisible 
for each $w$ dividing $p$.

\begin{theorem} \label{t1}
Assume that $H^{0}(F,\Ab(\rhob,r)) = H^{0}(F,\Ab(\rhob,r)^{*})=0$.
If $\mu(\rho_{0},r)=0$ for some $\rho_{0} \in \H(r)$, then
$\mu(\rho,r)=0$ for all $\rho \in \H(r)$.
\end{theorem}

We then say simply that $\mu(\rhob,r)=0$.

\begin{proof}
This is immediate from Corollary~\ref{cor:res}.
\end{proof}

\begin{theorem} \label{t2}
Assume that $H^{0}(F,\Ab(\rhob,r))=H^{0}(F,\Ab(\rhob,r)^{*})=0$ and
$\mu(\rhob,r)=0$.  The quantity
$$\lambda(\rho,r) - \sum_{w \nmid p} \dim_{k} \Ar^{\Gw}/\pi$$
is independent of the choice of $\rho \in \H(r)$.
\end{theorem}
\begin{proof}
This is immediate from Corollaries~\ref{cor:res} and~\ref{cor:quant}.
\end{proof}

Note that for a place $w$ of $\Fi$ dividing a place $v \nmid p$ of $F$, 
by \cite[Section 2]{GreenbergVatsal} the dimension
$\dim_{k} \Ar^{\Gw}/\pi$ is equal to the difference between the
number of occurrences of the Teichm\"uller character in the
$\Gw/\Iw$-representations $\Delta(\rho,r)_{\Iw}$ and
$A_{\rho,\Iw}[\pi]$.  Theorem~\ref{mt1} is thus an immediate
consequence of Theorems~\ref{t1} and~\ref{t2}.

\end{document}